\documentclass[12pt]{amsproc}
\usepackage[cp1251]{inputenc}
\usepackage[T2A]{fontenc} 
\usepackage[russian]{babel}
\usepackage{amsbsy,amsmath,amssymb,amsthm,graphicx,color} 
\usepackage{amscd}

\newtheorem{prop}{Предложение}

\theoremstyle{definition}

\theoremstyle{remark}

\begin {document}
\centerline{УДК 515.162.8}
\unitlength=1mm
\title[Функции роста групп Кокстера]
{Функции роста групп Кокстера, ряды Пуанкаре особенностей и $q$-дроби $2$-танглов}
\author{Г. Г. Ильюта}
\email{ilyuta@mccme.ru}
\address{}
\thanks{Работа поддержана грантами РФФИ-13-01-00755 и НШ-5138.2014.1}
\maketitle
%\begin{abstract}
%\end{abstract}

\bigskip

  Функции роста групп Кокстера, действующих отражениями в геодезических линиях в двумерных пространствах постоянной кривизны, представимы как частное многочленов Александера связной суммы $(2,k)$-торических узлов и крендельного узла \cite{5}. Рассматривая стандартные проекции этих узлов, замечаем, что они являются замыканиями одного и того же $2$-тангла -- суммы $2$-танглов, каждый из которых является степенью образующей группы кос из двух нитей. В статье получены близкие факты для функций роста всех групп Кокстера, для рядов Пуанкаре клейновых и фуксовых особенностей и для некоторых рядов Пуанкаре конечных групп. 

\begin{picture}(120,18)
\put(27,8){\circle{20}}
\put(57,8){\circle{20}}
\put(87,8){\circle{20}}

\put(29,10){\vector(1,1){4}}
\put(25,6){\vector(-1,-1){4}}
\put(33,2){\vector(-1,1){4}}
\put(21,14){\vector(1,-1){4}}

\put(59,10){\vector(1,1){4}}
\put(55,6){\vector(-1,-1){4}}
\put(63,2){\vector(-1,1){4}}
\put(51,14){\vector(1,-1){4}}

\put(89,10){\vector(1,1){4}}
\put(85,6){\vector(-1,-1){4}}
\put(93,2){\vector(-1,1){4}}
\put(81,14){\vector(1,-1){4}}

\qbezier(21,14)(27,19)(33,14)
\qbezier(21,2)(27,-3)(33,2)
\qbezier(81,14)(76,8)(81,2)
\qbezier(93,14)(98,8)(93,2)
\end{picture}

   Определим $2$-тангл как вложение двух отрезков и нескольких окружностей в трёхмерный шар, причём на границу попадают только концы отрезков. Замыкания $2$-тангла приводят к паре узлов, частное многочленов Александера которых назовём $q$-дробью $2$-тангла (многочлены нормализованы так, что замена переменной $q\to q^{-1}$ умножает многочлен на $\pm 1$ \cite{2}). 

\begin{prop}\label{prop1} Ряды Пуанкаре клейновых особенностей $E_6$, $E_7$, $E_8$ и $D_n$ представимы как $q$-дроби $2$-танглов.
\end{prop} 

  Фуксовы особенности связаны с действием фуксовых групп на касательном расслоении верхней полуплоскости. Если фуксова группа имеет сигнатуру $(g; a_1,\dots , a_r)$, то для ряда Пуанкаре $P(q)$ соответствующей особенности имеет место формула \cite{4}
$$
P(q)=\frac{1+(g-2)q+(g-2)q^2+q^3}{(1-q)^2}+\sum\nolimits_{i=1}^r\frac{q^2(1-q^{a_i-1})}{(1-q)^2(1-q^{a_i})}.
$$

\begin{prop}\label{prop2} Существует $2$-тангл, $q$-дробь которого равна $(1+q)^2P(-q)/q^{3/2}$.
\end{prop}

  Группа Кокстера $W_S$ порождается множеством инволюций $S=\{s_1,\dots , s_r\}$ и имеет определяющие соотношениями $(s_is_j)^{m_{ij}}=1$. Функция роста $W_S(q)$ группы Кокстера $W_S$ определяется как ряд 
$$
W_S(q)=\sum\nolimits_{g\in W_S}q^{l(g)},\quad l(g)=\min\{k: g=s_{i_1}\dots s_{i_k}\}.  
$$ 
Обозначим через $F_+(q)$ ($F_-(q)$) любой целочисленный многочлен, который делится на все многочлены $q^{\deg W_P(q)}+1$ ($q^{\deg W_P(q)}-1$), $P\subset S$, $1<|W_P|<\infty$, и частные являются возвратными многочленами. Пусть $z=q^{-1/2}-q^{1/2}$ и $\epsilon=0$ или $1$ в зависимости от того нечётна или чётна степень многочлена $F_{\pm}(q)$.

\begin{prop}\label{prop3} Для группы Кокстера $W_S$ функции
$$
\frac{z^\epsilon q^{\deg F_+/2}}{F_+(-q)} \left(\frac{1}{W_S(-q^{-1})}+\frac{1}{W_S(-q)}\right),   
$$
$$
\frac{z^\epsilon q^{\deg F_-/2}}{F_-(-q)} \left(\frac{1}{W_S(-q^{-1})}-\frac{1}{W_S(-q)}\right)   
$$
представимы как $q$-дроби $2$-танглов.
\end{prop}

  Для очень общего класса групп Кокстера функция роста удовлетворяет одному из равенств $W_S(q^{-1})=\pm W_S(q)$, в частности, если группа действует собственно и кокомпактно на стягиваемом $n$-многообразии (порождающие действуют как отражения), то $W_S(q^{-1})=(-1)^nW_S(q)$ \cite{1}.

  Эйлерова характеристика группы Кокстера обратна значению функции роста этой группы в $1$ \cite{10}. Поэтому доказанные факты связывают эйлерову характеристику и значение $q$-дроби в $-1$. Это значение известно как дробь $2$-тангла; обзор приложений и методов вычисления дроби $2$-тангла можно найти в \cite{8}. Модуль значения многочлена Александера в $-1$ известен как детерминант узла и он равен порядку группы гомологий двулистного накрытия сферы $S^3$ с ветвлением на узле. 

  Для конечной матричной группы $G\subset SU_n$ и её характера $\chi$ ряд Пуанкаре $G_\chi(q)$ можно определить с помощью формулы Молина
$$
G_\chi(q)=\sum\nolimits_{g\in G}\frac{\chi(g)}{\det(1-qg)}.   
$$

\begin{prop}\label{prop4} Если характер $\chi$ принимает целые значения и многочлены $\det(1-qg)$ являются целочисленными, то функции $q^{n/2}G_\chi(-q)/z$ для чётного $n$ и $q^{(n-1)/2}(1+q)G_\chi(-q)/z$ для нечётного $n$ представимы как $q$-дроби $2$-танглов.
\end{prop} 
 
  Доказательства. В \cite{3} ряды Пуанкаре клейновых особенностей представлены как частные многочленов Кокстера евклидовых и аффинных диаграмм Дынкина. Эти диаграммы в случаях $E_n$ и $D_n$ являются деревьями и их многочлены Кокстера реализуются как многочлены Александера слалом-узлов А'Кампо \cite{6}. Евклидова диаграмма Дынкина получается из аффинной удалением вершины. Узлы А'Кампо разложимы в композицию плюмбингов Хопфа, отвечающих вершинам соответствующих диаграмм Дынкина \cite{7} (для слалом-узлов наглядная интерпретация разложения на плюмбинги имеется в \cite{5}). Поэтому достаточно показать, что узлы, отличающиеся одним плюмбингом Хопфа, являются замыканиями одного и того же $2$-тангла. При плюмбинге Хопфа к поверхности Зейферта узла приклеивается лента. Разрезая эту ленту, мы получим $2$-тангл как границу поверхности (края разреза исключаются). Одно из замыканий этого $2$-тангла возвращает к исходному узлу (края разреза добавляются), а другое -- к узлу до разреза ленты (края разреза склеиваются). Процесс разложения слалом-узла в композицию плюмбингов Хопфа соответствует разложению частного многочленов Александера узла до и узла после одного плюмбинга в ветвящуюся цепную дробь \cite{6}. Эти процессы связаны теоремой Конвея: $q$-дробь суммы $2$-танглов равна сумме $q$-дробей слагаемых \cite{2}.

\begin{picture}(120,14)
\put(34,5){$+$}
\put(58,5){$=$}

\put(24,6){\circle{20}}
\put(48,6){\circle{20}}
\put(72,6){\circle{20}}
\put(96,6){\circle{20}}

\put(26,8){\vector(1,1){4}}
\put(22,4){\vector(-1,-1){4}}
\put(30,0){\vector(-1,1){4}}
\put(18,12){\vector(1,-1){4}}

\put(50,8){\vector(1,1){4}}
\put(46,4){\vector(-1,-1){4}}
\put(54,0){\vector(-1,1){4}}
\put(42,12){\vector(1,-1){4}}

\put(74,8){\vector(1,1){4}}
\put(70,4){\vector(-1,-1){4}}
\put(78,0){\vector(-1,1){4}}
\put(66,12){\vector(1,-1){4}}

\put(98,8){\vector(1,1){4}}
\put(94,4){\vector(-1,-1){4}}
\put(102,0){\vector(-1,1){4}}
\put(90,12){\vector(1,-1){4}}

\put(78,0){\line(1,0){12}}
\put(78,12){\line(1,0){12}}

\end{picture}

  В частности, $q$-дробью $2$-тангла является частное многочленов Александера $(2,n)$- и $(2,n+1)$-торических узлов, которые являются слалом-узлами диаграмм Дынкина $A_{n-1}$ и $A_n$. Поэтому для доказательства Предложения 2 остаётся представить как $q$-дробь $2$-тангла первое слагаемое в сумме
$$
\frac{(1+q)^2P(-q)}{q^{3/2}}=z(z^2-g+5)+\sum\nolimits_{i=1}^r\frac{q^{-(a_i-1)/2}-q^{(a_i-1)/2}}{q^{-a_i/2}-q^{a_i/2}}
$$
и применить теорему Конвея. В \cite{9} для любого целочисленного многочлена $b(z)$ построены $2$-танглы (один из другого получается отражением относительно биссектрисы первой четверти), $q$-дроби которых равны $zb(z^2)$ и $1/zb(z^2)$. Этот факт вместе с теоремой Конвея доказывает также Предложения 3 и 4. Функции из Предложения 3 представляются в виде суммы функций вида $1/zb(z^2)$ с помощью формулы Стейнберга \cite{12}
$$
\frac {1}{W_S(q^{-1})}=\sum\nolimits_{P\subset S,|W_P|<\infty}\frac {(-1)^{|P|}}{W_P(q)},   
$$
а функции из Предложения 4 -- с помощью формулы Молина. Это вытекает из следующих замечаний. Если группа $W_S$ конечна и её степени базисных инвариантов равны $m_1,\dots, m_r$, то $W_S(q)=\prod [m_i]$, где $[m]=1+q+\dots +q^{m-1}$ \cite{11}, в частности, многочлен $W_S(q)$ возвратный и делится на $1+q$, поскольку конечная группа Кокстера имеет инвариант степени $2$. Многочлен $F_+(q)$ делится на $1+q$, поскольку $\deg W_P(q)=1$, если $|P|=1$. Возвратный многочлен $R(q)$ нечётной (чётной) степени представляется в виде $q^{\deg R/2}R(-q)=zb(z^2)$ ($=b(z^2)$). Возвратный многочлен нечётной степени делится на $1+q$. Если возвратный многочлен чётной степени делится на $1+q$, то он делится на $(1+q)^2$. Вещественный характеристический многочлен $f(q)$ матрицы из $SU_n$ удовлетворяет равенству $q^nf(q^{-1})=(-1)^nf(q)$ и для нечётного $n$ делится на $1-q$. Многочлены $f(q)$ для чётного $n$ и $f(q)/(1-q)$ для нечётного $n$ являются возвратными.

  Приведём пример: для ряда Пуанкаре кольца инвариантов ($\chi\equiv 1$) бинарной группы тетраэдра формула Молина и соответствующее представление в виде суммы $q$-дробей из \cite{9} имеют вид
$$
G_1(q)=\frac{1}{(1-q)^2}+\frac{1}{(1+q)^2}+\frac{1}{1+q^2}+\frac{2}{1-q+q^2}+\frac{2}{1+q+q^2},   
$$
$$
\frac{qG_1(-q)}{z}=\frac{1}{z(z^2+4)}+\frac{1}{zz^2}+\frac{1}{z(z^2+2)}+2\frac{1}{z(z^2+3)}+2\frac{1}{z(z^2+1)}   
$$
и мы получаем другое представление ряда Пуанкаре особенности $E_6$ как $q$-дроби $2$-тангла.

\bigskip

\end {document}